\documentclass[12pt,leqno]{article}
\usepackage{hyperref}
\usepackage{soul}
\usepackage[doc]{optional}
\usepackage{xcolor}
\definecolor{labelkey}{rgb}{0,0.08,0.45}
\definecolor{rekey}{rgb}{0,0.6,0.0}
\definecolor{Brown}{rgb}{0.45,0.0,0.05}
\usepackage{exscale,relsize}
\usepackage{amsmath}
\usepackage{amsfonts}
\usepackage{amssymb}
\usepackage{calc}
\usepackage{theorem}
\usepackage{pifont}      
\usepackage{graphicx}
\newcommand{\scal}[2]{\langle{{#1},{#2}}\rangle}

\newcommand{\RX}{\ensuremath{\,\left]-\infty,+\infty\right]}}

\newcommand{\menge}[2]{\big\{{#1} \mid {#2}\big\}}

\newcommand{\To}{\ensuremath{\rightrightarrows}}

\newcommand{\dom}{\ensuremath{\operatorname{dom}}}

\newcommand{\gra}{\ensuremath{\operatorname{gra}}}

\newcommand{\ran}{\ensuremath{\operatorname{ran}}}

\newcommand{\pinf}{\ensuremath{+\infty}}

\renewcommand{\phi}{\ensuremath{\varphi}}

\newtheorem{theorem}{Theorem}[section]

\newtheorem{fact}[theorem]{Fact}
\newtheorem{corollary}[theorem]{Corollary}

\newtheorem{definition}[theorem]{Definition}

\theoremstyle{plain}{\theorembodyfont{\rmfamily}
}
\theoremstyle{plain}{\theorembodyfont{\rmfamily}
}
\theoremstyle{plain}{\theorembodyfont{\rmfamily}
}
\theoremstyle{plain}{\theorembodyfont{\rmfamily}
}
\theoremstyle{plain}{\theorembodyfont{\rmfamily}
\newtheorem{remark}[theorem]{Remark}}

\theoremstyle{plain}{\theorembodyfont{\rmfamily}
}

\begin{document}


\title{\sffamily{Every maximally monotone operator of \\
Fitzpatrick-Phelps type is actually of dense type}}

\author{
Heinz H.\ Bauschke\thanks{Mathematics, Irving K.\ Barber School,
University of British Columbia, Kelowna, B.C. V1V 1V7, Canada. E-mail:
\texttt{heinz.bauschke@ubc.ca}.},\;
Jonathan M. Borwein\thanks{CARMA, University of Newcastle, Newcastle, New South Wales 2308, Australia. E-mail:
\texttt{jonathan.borwein@newcastle.edu.au}.},\;
 Xianfu
Wang\thanks{Mathematics, Irving K.\ Barber School, University of British Columbia,
Kelowna, B.C. V1V 1V7, Canada. E-mail:
\texttt{shawn.wang@ubc.ca}.},\; and Liangjin\
Yao\thanks{Mathematics, Irving K.\ Barber School, University of British Columbia,
Kelowna, B.C. V1V 1V7, Canada.
E-mail:  \texttt{ljinyao@interchange.ubc.ca}.}}

\date{April 4, 2011}
\maketitle

\begin{abstract} \noindent
We show that every maximally monotone operator of Fitzpatrick-Phelps
type defined on a real Banach space must be of dense type. This
provides an affirmative answer to a question posed by Stephen Simons
in 2001 and implies that various important notions of monotonicity
coincide.
\end{abstract}

\noindent {\bfseries 2010 Mathematics Subject Classification:}\\
{Primary   47H05;
Secondary
46B10, 47N10,
 90C25}

\noindent {\bfseries Keywords:}
Fitzpatrick function,
maximally monotone operator,
monotone operator,
multifunction,
operator of type (D),
operator of type (FP),
operator of type (NI),
set-valued operator.

\section{Introduction}

Throughout this note, we assume that
$X$ is a real Banach space with norm $\|\cdot\|$,
that $X^*$ is the continuous dual of $X$, and
that $X$ and $X^*$ are paired by $\scal{\cdot}{\cdot}$.
Let $A\colon X\To X^*$
be a \emph{set-valued operator} (also known as multifunction)
from $X$ to $X^*$, i.e., for every $x\in X$, $Ax\subseteq X^*$,
and let
$\gra A = \menge{(x,x^*)\in X\times X^*}{x^*\in Ax}$ denote
the \emph{graph} of $A$. The \emph{domain} of $A$ is
$\dom A= \menge{x\in X}{Ax\neq\varnothing}$, while
$\ran A=A(X)$ is the \emph{range} of $A$.
Recall that $A$ is  \emph{monotone} if
\begin{equation}
\scal{x-y}{x^*-y^*}\geq 0,\quad \forall (x,x^*)\in \gra A\;\forall (y,y^*)\in\gra A,
\end{equation}
and \emph{maximally monotone} if $A$ is monotone and $A$ admits
 no proper monotone extension.
It will be convenient to also say that $\gra A$ is monotone or
maximally monotone respectively in this case.
We can then simply say that
$(x,x^*)\in X\times X^*$ is monotonically related
to $\gra A$ if $\{(x,x^*)\}\cup\gra A$ is monotone.

 We now recall the three fundamental types of monotonicity.

 \begin{definition}\label{def1}
 Let $A:X\To X^*$ be maximally monotone.
 Then three key types of monotone operators are defined as follows.
 \begin{enumerate}
 \item $A$ is
\emph{of dense type or type (D)} (1971, \cite{Gossez3}) if for every
$(x^{**},x^*)\in X^{**}\times X^*$ with
\begin{align*}
\inf_{(a,a^*)\in\gra A}\langle a-x^{**}, a^*-x^*\rangle\geq 0,
\end{align*}
there exist a  bounded net
$(a_{\alpha}, a^*_{\alpha})_{\alpha\in\Gamma}$ in $\gra A$
such that
$(a_{\alpha}, a^*_{\alpha})_{\alpha\in\Gamma}$
weak*$\times$strong converges to
$(x^{**},x^*)$.
\item $A$ is
\emph{of type negative infimum (NI)} (1996, \cite{SiNI}) if
\begin{align*}
\sup_{(a,a^*)\in\gra A}\big(\langle a,x^*\rangle+\langle a^*,x^{**}\rangle
-\langle a,a^*\rangle\big)
\geq\langle x^{**},x^*\rangle,
\quad \forall(x^{**},x^*)\in X^{**}\times X^*.
\end{align*}
\item
$A$ is \emph{of type Fitzpatrick-Phelps (FP)} (1992, \cite{FP}) if
whenever $U$ is an open convex subset of $X^*$ such that $U\cap \ran
A\neq\varnothing$, $x^*\in U$, and $(x,x^*)\in X\times X^*$ is
monotonically related to $\gra A\cap (X\times U)$ it must follow
that $(x,x^*)\in\gra A$.
\end{enumerate}
\end{definition}

All three of these properties are known to hold for the subgradient
of a closed convex function and for every maximally monotone
operator on a reflexive space. These and other relationships known
amongst these and other monotonicity notions are described in
\cite[Chapter 8]{BorVan}. Monotone operators are fundamental objects
in modern Optimization and Variational Analysis; see, e.g.,
\cite{Bor1,Bor2,Bor3}, the books \cite{BC2011,
BorVan,BurIus,ph,Si,Si2,RockWets,Zalinescu} and the references
therein.

In Theorem \ref{FP:1} of  this paper, we provide an affirmative to
the following question, posed by S.\ Simons
\cite[Problem~18,~page~406]{Si4}:

\begin{quote} \emph{Let $A:X\To X^*$ be maximally monotone such that $A$ is
of type (FP).\\  Is $A$ necessarily of type (D)?}\end{quote}

In consequence, in Corollary \ref{cor:main} we record that the three
notions in Definition \ref{def1} actually coincide.

We shall utilize the following notation, in addition to standard
notions from convex analysis:
The \emph{open unit ball} in $X$ is $U_X=
\menge{x\in X}{\|x\|<1}$, and the \emph{closed unit ball}
is $\menge{x\in X}{\|x\|\leq 1}$.
It is very convenient to identify $X$ with its canonical image in
the bidual space $X^{**}$. Moreover, $X\times X^*$ and $(X\times
X^*)^* = X^*\times X^{**}$ are paired via
$$\scal{(x,x^*)}{(y^*,y^{**})} = \scal{x}{y^*} +
\scal{x^*}{y^{**}},$$ where $(x,x^*)\in X\times X^*$ and
$(y^*,y^{**}) \in X^*\times X^{**}$. We recall the following basic
fact regarding the second dual ball:

\begin{fact}[Goldstine]\emph{(See \cite[Theorem~2.6.26]{Megg} or \cite[Theorem~3.27]{FabianHH}.)}
  \label{Goldst:1}
 The weak*-closure of $B_X$ in $X^{**}$ is $B_{X^{**}}$.
\end{fact}

The remainder of this paper is organized as follows. In
Section~\ref{s:aux}, we record auxiliary results for subsequent use.
The main result (Theorem~\ref{FP:1}) and the promised corollary
(Corollary~\ref{cor:main}) are provided in Section~\ref{s:main}.

 \section{Preliminary monotonicity results}\label{s:aux}

A now fundamental tool of modern monotone operator theory originated
with Simon Fitzpatrick  in 1988. It is reprised next:

\begin{fact}[Fitzpatrick]
\emph{(See {\cite[Corollary~3.9]{Fitz88}}.)}
\label{f:Fitz}
Let $A\colon X\To X^*$ be maximally monotone,  and let us set
\begin{equation}
F_A\colon X^{**}\times X^*\to\RX\colon
(x^{**},x^*)\mapsto \sup_{(a,a^*)\in\gra A}
\big(\scal{x^{**}}{a^*}+\scal{a}{x^*}-\scal{a}{a^*}\big).
\end{equation}
Then for every $(x,x^*)\in X\times X^*$, the inequality
$$\scal{x}{x^*}\leq F_A(x,x^*)$$ is true,
and equality holds if and only if $(x,x^*)\in\gra A$.\\ The function
${F_A}|_{X\times X^*}$ is the classical \emph{Fitzpatrick function}
associated with $A$.
\end{fact}

The first relevant relationship established for (FP) operators is
due to Stephen Simons:

\begin{fact}[Simons]
\emph{(See \cite[Theorem~17]{Si4} or \cite[Theorem~37.1]{Si2}.)}
\label{FTSim:1}
Let $A:X\To X^*$ be maximally monotone
and of type (D). Then $A$ is of type (FP).
\end{fact}

The most powerful current information is captured in the following
result.

\begin{fact}[Simons / Marques Alves and Svaiter]\emph{(See
\cite[Lemma~15]{SiNI} or
\cite[Theorem~36.3(a)]{Si2}, and \cite[Theorem~4.4]{MarSva}.)}
\label{PF:Su1} Let $A:X \To X^*$ be maximally  monotone.
Then $A$ is of type (D) if and only if it is of type (NI).
\end{fact}

The implication type (NI) implies type (D) --- which we exploit
--- is very recently due to Marques Alves and Svaiter \cite{MarSva}.

\section{Main result}\label{s:main}

The next theorem is our main result. In conjunction with the
corollary that follows, it provides the affirmative answer  promised
to Simons's problem posed in \cite[Problem~18]{Si4}.

\begin{theorem}\label{FP:1}
Let $A:X\To X^*$ be maximally monotone such that
$A$ is of type (FP).
Then $A$ is of type (NI).
\end{theorem}
\begin{proof}
After translating the graph if necessary,
we can and do suppose that $(0,0)\in\gra A$.
Let $(x_0^{**},x^*_0)\in X^{**}\times X^*$.
We must show that
\begin{align}
F_A(x_0^{**},x^*_0)\geq \langle x_0^{**},x_0^*\rangle\label{FPSe:1}
\end{align}
and we consider two cases.

\emph{Case 1}: $x_0^{**}\in X$.\\
Then \eqref{FPSe:1} follows  directly from Fact~\ref{f:Fitz}.

\emph{Case 2}: $x_0^{**}\in X^{**}\smallsetminus X$.\\
By Fact~\ref{Goldst:1},  there exists a bounded net $(x_\alpha)_{\alpha\in I}$
in $X$  that weak* converges to $x_0^{**}$.
Thus, we have
\begin{align}
M = \sup_{\alpha \in I} \|x_\alpha\|<\pinf\label{FPSe:1a}
\end{align}
and
\begin{align}
\langle x_\alpha, x^*_0\rangle\rightarrow\langle x_0^{**},x^*_0\rangle.\label{FPSe:1aa}
\end{align}
Now we consider two subcases.

 \emph{Subcase 2.1}: There exists $\alpha\in I$,
  such that
 $(x_{\alpha},x^*_0)\in\gra A$.\\
By definition,
\begin{align*}
&F_A(x_0^{**},x_0^*)
\geq \langle x_{\alpha},x_0^*\rangle+\langle x_0^{**},x_0^*\rangle
-\langle x_{\alpha},x_0^*\rangle
=\langle x_0^{**},x_0^*\rangle.
\end{align*}
Hence \eqref{FPSe:1} holds.

 \emph{Subcase 2.2}:
 We have
\begin{align}
(x_{\alpha},x^*_0)\notin\gra A,\quad \forall \alpha\in I.\label{FPSe:2}
\end{align}
Set \begin{align}
U_{\varepsilon}=\left[0, x^*_0\right]+\varepsilon U_{X^*},\label{FPSe:3}
\end{align} where $\varepsilon>0$.
Observe that $U_{\varepsilon}$ is open and convex.
Since $(0,0)\in\gra A$, we have, by definition of $U_\varepsilon$,
$0\in\ran A\cap U_{\varepsilon}$ and $x^*_0\in U_{\varepsilon}$.
In view of \eqref{FPSe:2} and because $A$ is of type (FP),
there exists a net
$(a_{\alpha,\varepsilon}, a^*_{\alpha,\varepsilon})$ in $\gra A$
such that $ a^*_{\alpha,\varepsilon}\in U_{\varepsilon}$ and
\begin{align}
\langle a_{\alpha,\varepsilon}, x^*_0\rangle+\langle x_{\alpha},  a^*_{\alpha,\varepsilon}\rangle
-\langle a_{\alpha,\varepsilon}, a^*_{\alpha,\varepsilon}\rangle>\langle x_{\alpha},x^*_0\rangle,\quad \forall \alpha\in I.\label{FPSe:4}
\end{align}
Now fix $\alpha\in I$.
By \eqref{FPSe:4},
\begin{align}
&\langle a_{\alpha,\varepsilon}, x^*_0\rangle+\langle x_0^{**},  a^*_{\alpha,\varepsilon}\rangle
-\langle a_{\alpha,\varepsilon}, a^*_{\alpha,\varepsilon}\rangle
>\langle x_0^{**}-x_{\alpha},  a^*_{\alpha,\varepsilon}\rangle+\langle x_{\alpha},x^*_0\rangle.\nonumber\end{align}
Hence,
\begin{align}
F_A(x_0^{**},x_0^*)>\langle x_0^{**}-x_{\alpha},
a^*_{\alpha,\varepsilon}\rangle+\langle x_{\alpha},x^*_0\rangle.
\label{FPSe:5}
\end{align}
Since $a^*_{\alpha,\varepsilon}\in U_{\varepsilon}$,
there exist
\begin{align}
t_{\alpha,\varepsilon}\in\left[0,1\right]\ \text{and}\ b^*_{\alpha,\varepsilon}\in U_{X^*}\label{FPSe:06}\end{align}
 such that
\begin{align}
a^*_{\alpha,\varepsilon}=t_{\alpha,\varepsilon}x^*_0+\varepsilon b^*_{\alpha,\varepsilon}.\label{FPSe:6}
\end{align}
Using \eqref{FPSe:5}, \eqref{FPSe:6},
and \eqref{FPSe:1a}, we deduce that
\begin{align}
F_A(x_0^{**},x_0^*)&>\langle x_0^{**}
-x_{\alpha},  t_{\alpha,\varepsilon}x^*_0+\varepsilon b^*_{\alpha,\varepsilon}\rangle+\langle x_{\alpha},x^*_0\rangle\nonumber\\
&=t_{\alpha,\varepsilon}\langle x_0^{**}-x_{\alpha},  x^*_0\rangle+\varepsilon
\langle x_0^{**}-x_{\alpha}, b^*_{\alpha,\varepsilon}\rangle+\langle x_{\alpha},x^*_0\rangle\nonumber\\
&\geq t_{\alpha,\varepsilon}\langle x_0^{**}-x_{\alpha},  x^*_0\rangle-\varepsilon
\| x_0^{**}-x_{\alpha}\|+\langle x_{\alpha},x^*_0\rangle\nonumber\\
&\geq t_{\alpha,\varepsilon}\langle x_0^{**}-x_{\alpha},  x^*_0\rangle-\varepsilon
(\| x_0^{**}\|+M)+\langle x_{\alpha},x^*_0\rangle.
\label{FPSe:7}
\end{align}
In view of \eqref{FPSe:06} and since $\alpha\in I$ was
chosen arbitrarily, we take  the limit in \eqref{FPSe:7}
and obtain with the help of \eqref{FPSe:1aa} that
\begin{align}
F_A(x_0^{**},x_0^*)\geq-\varepsilon
(\| x_0^{**}\|+M)+\langle x_0^{**},x^*_0\rangle.\label{FPSe:8}
\end{align}
Next, letting $\varepsilon\rightarrow 0$ in \eqref{FPSe:8}, we have
 \begin{align}
F_A(x_0^{**},x_0^*)\geq\langle x_0^{**},x^*_0\rangle.\label{FPSe:9}
\end{align}
Therefore,
\eqref{FPSe:1} holds in all cases.
\end{proof}

We now obtain the promised corollary:

\begin{corollary}
\label{cor:main} Let $A\colon X\To X^*$ be maximally monotone. Then
the following are equivalent.
\begin{enumerate}
\item\label{MaFT:1} $A$ is of type (D).
\item \label{MaFT:2} $A$ is of type (NI).
\item \label{MaFT:3} $A$ is of type (FP).
\end{enumerate}
\end{corollary}
\begin{proof}
First \emph{(i) implies (iii)} is Fact~\ref{FTSim:1}; next
Theorem~\ref{FP:1}  shows \emph{(iii) implies (ii)}; while
Fact~\ref{PF:Su1} implies  concludes the circle with \emph{(ii)
implies (i)}.
\end{proof}

We note that while the result is now quite easy, it remained
inaccessible until \cite[Theorem~4.4]{MarSva} was available.

\begin{remark}
Let $A\colon X\To X^*$ be maximally monotone.
Corollary~\ref{cor:main} establishes the equivalences of the key
types (D), (NI), and (FP), which as noted all hold when $X$ is
reflexive or $A=\partial f$, where $f\colon X\to\RX$ is convex,
lower semicontinuous, and proper (see \cite{BorVan,Si,Si2}).

Furthermore, these notions are also equivalent to type (ED), see
\cite{Si5}. For a nonlinear operator they also coincide with
\emph{uniqueness} of maximal extensions to $X^{**}$ (see
\cite{MarSva}). In \cite[p.~454]{BorVan} there is discussion of this
result and of the linear case.

Finally, when $A$ is a linear relation, it has recently been
established that all these notions coincide with monotonicity of the
\emph{adjoint} multifunction  $A^*$ (see \cite{BBWY:1}).
\end{remark}

\section*{Acknowledgments}
Heinz Bauschke was partially supported by the Natural Sciences and
Engineering Research Council of Canada and
by the Canada Research Chair Program.
Jonathan  Borwein was partially supported by the Australian Research  Council.
Xianfu Wang was partially supported by the Natural
Sciences and Engineering Research Council of Canada.

\end{document}